\documentclass[11pt,final]{amsart}
\usepackage{amsmath,amssymb,amsthm,showkeys}

\usepackage{amsmath,amssymb,amsthm,showkeys}
\newif\ifpdf
    \ifx\pdfoutput\undefined
    \pdffalse 
    \else
    \pdfoutput=1 
    \pdftrue
    \fi

    \ifpdf
    \usepackage[pdftex]{graphicx}
    \else
    \usepackage{graphicx}
    \fi

\usepackage{amsmath}
\usepackage{graphicx}
\usepackage{color}

\begin{document}
\ifpdf
    \DeclareGraphicsExtensions{.pdf, .jpg}
    \else
    \DeclareGraphicsExtensions{.eps, .jpg}
    \fi

\newtheorem{theorem}{Theorem}
 \newtheorem{lemma}{Lemma}
 \newtheorem{proposition}{Proposition}
 \newtheorem{remark}{Remark}
 \newtheorem{question}{Question}
 \newtheorem{definition}{Definition}
 \newtheorem{corollary}{Corollary}
 \newtheorem{claim}{Claim}
 \newtheorem{notation}{Notation}
\newtheorem{lem}{Lemma}[section]
\newtheorem{prop}{Proposition}[section]
\newtheorem{cor}{Corollary}[section]
\numberwithin{equation}{section}
\newtheorem{thm}{Theorem}[section]

\newtheorem{ass}{Assumption}[section]

\theoremstyle{remark}
\newtheorem{example}{Example}[section]
\newtheorem*{ack}{Acknowledgments}

\theoremstyle{definition}

\theoremstyle{remark}
\theoremstyle{remark}

\newenvironment{Abstract}
{\begin{center}\textbf{\footnotesize{Abstract}}%
\end{center} \begin{quote}\begin{footnotesize}}
{\end{footnotesize}\end{quote}\bigskip}

\newcommand{\triple}[1]{{|\!|\!|#1|\!|\!|}}

\newcommand{\xx}{\langle x\rangle}
\newcommand{\ep}{\varepsilon}
\newcommand{\al}{\alpha}
\newcommand{\be}{\beta}
\newcommand{\de}{\partial}
\newcommand{\la}{\lambda}
\newcommand{\La}{\Lambda}
\newcommand{\ga}{\gamma}
\newcommand{\del}{\delta}
\newcommand{\Del}{\Delta}
\newcommand{\sig}{\sigma}
\newcommand{\ome}{\omega}
\newcommand{\Ome}{\Omega}
\newcommand{\C}{{\mathbb C}}
\newcommand{\N}{{\mathbb N}}
\newcommand{\Z}{{\mathbb Z}}
\newcommand{\R}{{\mathbb R}}
\newcommand{\Rn}{{\mathbb R}^{n}}
\newcommand{\Rnu}{{\mathbb R}^{n+1}_{+}}
\newcommand{\Cn}{{\mathbb C}^{n}}
\newcommand{\spt}{\,\mathrm{supp}\,}
\newcommand{\Lin}{\mathcal{L}}
\newcommand{\SSS}{\mathcal{S}}
\newcommand{\F}{\mathcal{F}}
\newcommand{\xxi}{\langle\xi\rangle}
\newcommand{\eei}{\langle\eta\rangle}
\newcommand{\xei}{\langle\xi-\eta\rangle}
\newcommand{\yy}{\langle y\rangle}
\newcommand{\dint}{\int\!\!\int}
\newcommand{\hatp}{\widehat\psi}

\newcommand{\q}{\quad}
\newcommand{\qq}{\quad\quad}
\newcommand{\qqq}{\quad\quad\quad}
\newcommand{\qqqq}{\quad\quad\quad\quad}
\newcommand{\qqqqq}{\quad\quad\quad\quad\quad}
\newcommand{\qqqqqq}{\quad\quad\quad\quad\quad\quad}

\newcommand{\abs}[1]{\lvert#1\rvert}
\newcommand{\norm}[2]{{\left\| #1 \right\|}_{#2}}
\newcommand{\f}[2]{\frac{#1}{#2}}
\newcommand{\dpr}[2]{\langle #1,#2 \rangle}

\newcommand{\ppr}{p'}
\newcommand{\pp}[2]{\f{\p #1}{\p #2}}
\newcommand{\pnf}{+\infty}
\newcommand{\mnf}{-\infty}
\newcommand{\nf}{\infty}
\newcommand{\oop}{\frac{1}{p}}
\newcommand{\ooppr}{\frac{1}{p'}}
\newcommand{\ooq}{\frac{1}{q}}
\newcommand{\ooqpr}{\frac{1}{q'}}
\newcommand{\oor}{\frac{1}{r}}
\newcommand{\oorpr}{\frac{1}{r'}}

\newcommand{\Ga}{\Gamma}
\newcommand{\ben}{\beta_n}
\newcommand{\De}{\Delta}
\newcommand{\ve}{\varepsilon}
\newcommand{\Ve}{\Xi}
\newcommand{\ze}{\zeta}
\newcommand{\Th}{\Theta}
\newcommand{\ka}{\kappa}
\newcommand{\laj}{\lambda_j}
\newcommand{\lak}{\lambda_k}
\newcommand{\si}{\sigma}
\newcommand{\Si}{\Sigma}
\newcommand{\vp}{\varphi}
\newcommand{\om}{\omega}
\newcommand{\Om}{\Omega}

\newcommand{\rn}{{\mathbf R}^n}
\newcommand{\rmm}{{\mathbf R}^m}
\newcommand{\rone}{\mathbf R^1}
\newcommand{\rtwo}{\mathbf R^2}
\newcommand{\rthree}{\mathbf R^3}
\newcommand{\rfour}{\mathbf R^4}
\newcommand{\ronen}{{\mathbf R}^{n+1}}

\newcommand{\tn}{\mathbf T^n}
\newcommand{\tone}{\mathbf T^1}
\newcommand{\ttwo}{\mathbf T^2}
\newcommand{\tthree}{\mathbf T^3}
\newcommand{\tfour}{\mathbf T^4}

\newcommand{\zn}{\mathbf Z^n}
\newcommand{\zp}{\mathbf Z^+}
\newcommand{\zone}{\mathbf Z^1}
\newcommand{\zz}{\mathbf Z}
\newcommand{\ztwo}{\mathbf Z^2}
\newcommand{\zthree}{\mathbf Z^3}
\newcommand{\zfour}{\mathbf Z^4}

\newcommand{\hn}{\mathbf H^n}
\newcommand{\hone}{\mathbf H^1}
\newcommand{\htwo}{\mathbf H^2}
\newcommand{\hthree}{\mathbf H^3}
\newcommand{\hfour}{\mathbf H^4}

\newcommand{\cone}{\mathbf C^1}
\newcommand{\ctwo}{\mathbf C^2}
\newcommand{\cthree}{\mathbf C^3}
\newcommand{\cfour}{\mathbf C^4}

\newcommand{\sn}{\mathbf S^{n-1}}
\newcommand{\sone}{\mathbf S^1}
\newcommand{\stwo}{\mathbf S^2}
\newcommand{\sthree}{\mathbf S^3}
\newcommand{\sfour}{\mathbf S^4}

\newcommand{\lp}{L^{p}}
\newcommand{\lppr}{L^{p'}}
\newcommand{\lqq}{L^{q}}
\newcommand{\lr}{L^{r}}

\newcommand{\wlp}{L^{p,\infty}}
\newcommand{\wlq}{L^{q,\infty}}
\newcommand{\wlr}{L^{r,\infty}}
\newcommand{\wlo}{L^{1,\infty}}

\newcommand{\lprn}{L^{p}(\rn)}
\newcommand{\lptn}{L^{p}(\tn)}
\newcommand{\lpzn}{L^{p}(\zn)}
\newcommand{\lpcn}{L^{p}(\cn)}
\newcommand{\lphn}{L^{p}(\cn)}

\newcommand{\lprone}{L^{p}(\rone)}
\newcommand{\lptone}{L^{p}(\tone)}
\newcommand{\lpzone}{L^{p}(\zone)}
\newcommand{\lpcone}{L^{p}(\cone)}
\newcommand{\lphone}{L^{p}(\hone)}

\newcommand{\lqrn}{L^{q}(\rn)}
\newcommand{\lqtn}{L^{q}(\tn)}
\newcommand{\lqzn}{L^{q}(\zn)}
\newcommand{\lqcn}{L^{q}(\cn)}
\newcommand{\lqhn}{L^{q}(\hn)}

\newcommand{\lo}{L^{1}}
\newcommand{\lt}{L^{2}}
\newcommand{\li}{L^{\infty}}

\newcommand{\lmixpq}{L^p_tL^q_x}

\newcommand{\co}{C^{1}}
\newcommand{\ci}{C^{\infty}}
\newcommand{\coi}{C_0^{\infty}}

\newcommand{\cs}{\mathcal S}
\newcommand{\cm}{\mathcal M}
\newcommand{\cf}{\mathcal F}
\newcommand{\cb}{\mathcal B}
\newcommand{\ce}{\mathcal E}
\newcommand{\cd}{\mathcal D}
\newcommand{\csp}{\mathcal S'}
\newcommand{\cn}{\mathcal N}
\newcommand{\cz}{\mathcal Z}
\newcommand{\crr}{\mathbf R}
\newcommand{\cc}{\mathcal C}
\newcommand{\ch}{\mathcal H}
\newcommand{\cl}{\mathcal L}

\newcommand{\pv}{\textup{p.v.}\,}
\newcommand{\loc}{\textup{loc}}
\newcommand{\intl}{\int\limits}
\newcommand{\iintl}{\iint\limits}
\newcommand{\diint}{\displaystyle\iint}
\newcommand{\dintl}{\displaystyle\intl}
\newcommand{\diintl}{\displaystyle\iintl}
\newcommand{\liml}{\lim\limits}
\newcommand{\suml}{\sum\limits}
\newcommand{\supl}{\sup\limits}
\newcommand{\df}{\displaystyle\frac}
\newcommand{\p}{\partial}
\newcommand{\Ar}{\textup{Arg}}
\newcommand{\abssigk}{\widehat{|\si_k|}}
\newcommand{\tT}{\tilde{T}}
\newcommand{\tV}{\tilde{V}}
\newcommand{\wt}{\widetilde}
\newcommand{\Qvi}{Q_{\nu,i}}
\newcommand{\sjv}{a_{j,\nu}}
\newcommand{\sja}{a_{j,1}}
\newcommand{\sj}{a_j}
\newcommand{\pvs}{P_\nu^s}
\newcommand{\pva}{P_1^s}
\newcommand{\cjk}{c_{j,k}^{m,s}}
\newcommand{\Bjsnu}{B_{j-s,\nu}}
\newcommand{\Bjs}{B_{j-s}}
\newcommand{\Ly}{L_i^y}
\newcommand{\dd}[1]{\f{\partial}{\partial #1}}
\newcommand{\czz}{Calder\'on-Zygmund }

\newcommand{\lbl}{\label}
\newcommand{\beq}{\begin{equation}}
\newcommand{\eeq}{\end{equation}}
\newcommand{\beqna}{\begin{eqnarray*}}
\newcommand{\eeqna}{\end{eqnarray*}}
\newcommand{\beqn}{\begin{equation*}}
\newcommand{\eeqn}{\end{equation*}}
\newcommand{\bp}{\begin{proof}}
\newcommand{\bprop}{\begin{proposition}}
\newcommand{\eprop}{\end{proposition}}
\newcommand{\bt}{\begin{theorem}}
\newcommand{\et}{\end{theorem}}
\newcommand{\bex}{\begin{Example}}
\newcommand{\eex}{\end{Example}}
\newcommand{\bc}{\begin{corollary}}
\newcommand{\ec}{\end{corollary}}
\newcommand{\bcl}{\begin{claim}}
\newcommand{\ecl}{\end{claim}}
\newcommand{\bl}{\begin{lemma}}
\newcommand{\el}{\end{lemma}}

\renewcommand{\Re}{\;\mathrm{Re}\;}
\renewcommand{\Im}{\;\mathrm{Im}\;}

\title[Smoothing - Strichartz estimates]{Smoothing -  Strichartz
Estimates for the  Schr\"odinger Equation with small Magnetic
Potential}

\author{Vladimir Georgiev, Atanas Stefanov and Mirko Tarulli}



\thanks{V.G. and M.T.  are
  partially
supported by Research Training Network (RTN) HYKE,  financed by
the European Union, contract number : HPRN-CT-2002-00282.
A.S. is supported in part
by an NSF grant  0300511.}

\maketitle

\date{}


\begin{Abstract}The work treats smoothing and dispersive
properties of solutions to the Schr\"odinger equation with
magnetic potential. Under suitable smallness assumption on the
potential involving scale invariant norms we prove smoothing -
Strichartz estimate for the corresponding Cauchy problem. An
application that guarantees absence of pure point spectrum of the
corresponding perturbed Laplace operator is discussed too.
\end{Abstract}

\vskip+0.2cm

 \vskip+0.5cm

\section{Introduction and statement of results}

Let $A = (A_1(t,x),\cdots,A_n(t,x)), x\in \mathbb{R}^n, n \geq 3$
 be a
magnetic potential, such that $A_j(t,x), j=1,\cdots,n,$ are real
valued functions, and let the magnetic Laplacian operator be
$$
\Delta_A=\sum_j  (\p_j+i A_j)^2  = \De+ 2 i A\nabla  + i div(A) - (\sum_j A_j^2)
$$
Our goal is to study the dispersive properties
of the corresponding Schr\"odinger equation
\begin{equation}\label{eq.lwa}
\begin{cases}
\partial_t u - i\Delta_A u = F(t,x) ,\ \ \  t \in \mathbb{R} , \ \ \ x \in \mathbb{R}^n\\
u(0,x)=f(x).
\end{cases}
\end{equation}
In this paper, we will be concerned with
the Strichartz and smoothing estimates for \eqref{eq.lwa}, when the vector potential $A$ is small in certain sense.
In fact, we aim at obtaining global scale invariant
Strichartz and smoothing estimates, under appropriate scale invariant smallness assumptions on $A$.

In the ``free'' case $A=0$, there exists vast literature on the subject.
Introduce the mixed space-time norms
$$
\| u \|_{L^q_t L^{r}_x} = \left(\int_{\mathbb{R}} \left(
\int_{\mathbb{R}^n} |u(t,x)|^r dx\right)^{q/r} dt \right)^{1/q}.
$$
We say that a pair of exponents $(q,r)$ is Strichartz admissible,
if $2\leq  q,r\leq\infty$, $2/q+n/r=n/2$ and $(q,r,n)\neq (2,\infty, 2)$.
Then, by  result of
Strichartz,
Ginibre-Velo,  and  Keel-Tao,
\begin{eqnarray}
\label{eq:200}
& & \|e^{ i t \De} f\|_{L^q_t L^r_x}\leq C\|f\|_{L^2}\\
\label{eq:201} & & \left\|\intl e^{ i s  \De}
F(s,\cdot) ds\right\|_{L^2_x}\leq C\norm{F}{L^{\tilde{q}'}_t L^{\tilde{r}'}_x} \\
& & \label{eq:202} \left\|\intl_0^t e^{ i (t-s) \De} F(s,\cdot)
ds\right\|_{L^q_t L^r_x}\leq C\norm{F}{L^{\tilde{q}'}_t
L^{\tilde{r}'}_x},
\end{eqnarray}
where $(\tilde{q}, \tilde{r})$ is another Strichartz admissible
pair and $q'=q/(q-1)$. Note that for $n \geq 3$ the set of
admissible pairs $(q, r)$ can be represented equivalently as
$(1/q,1/r) \in AB,$ where $AB$ is the segment with end points
$A(0,1/2),$ $B(1/2, 2n/(n-2))$ and we can rewrite the estimate
\eqref{eq:201} as
\begin{equation}\label{eq.201eq}
    \left\|\intl e^{ i s  \De}
F(s,\cdot) ds\right\|_{L^2_x}\leq C \left( \inf_{F=F_1+F_2}
\norm{F_1}{L^1_t L^2_x} + \norm{F_2}{L^2_t L^{2n/(n+2)}_x}
\right).
\end{equation}

On the other hand, the smoothing estimates were established by
Kenig-Ponce-Vega in the  seminal paper, \cite{KPV93}, see also
Ruiz-Vega \cite{RuVe94} . These were later extended to more
general second order Schr\"odinger equations in \cite{KPV98}. Some
possible scale and rotation invariant smoothing estimates similar
to \eqref{eq:200}, \eqref{eq:201} and \eqref{eq:202} can be
written as (see Corollary \ref{corsmvs}  below   )
\begin{equation}\label{eq.3214}
   \supl_{m \in \Z} \ \left(2^{-m/2} 2^{k/2} \norm{ e^{ i t \De} f_k}{L^2_t L^2(|x|\sim
2^{m})}\right) \leq C\norm{f_k}{L^2},
\end{equation}
\begin{equation}\label{eq.sm2323}
    \left\|\intl e^{ i s  \De}
F_k(s,\cdot) ds\right\|_{L^2_x}\leq C \left( \sum_{m \in \Z}
2^{m/2} 2^{-k/2} \norm{  F_k}{L^2_t L^2(|x|\sim 2^{m})} \right).
\end{equation}
\begin{eqnarray}\label{eq.sminh2323}
   \supl_{m \in \Z} \left( 2^{-m/2} 2^{k/2} \left\| \intl_0^t e^{ i (t-s) \De} F_k(s,\cdot) ds\right\|_{L^2_t L^2(|x|\sim
   2^{m})}\right) \
   \leq \\ \nonumber
   \leq  C \left( \sum_{m \in \Z} 2^{m/2}2^{-k/2} \norm{
 F_k}{L^2_t L^2(|x|\sim 2^{m})} \right),
\end{eqnarray}
where $k$ is any integer,  $\phi_k:=P_k \phi$ is the $k^{th}$
Littlewood-Paley piece of $\phi$ ( see Section \ref{sec:2.1}
below).

Motivated by  these estimates, given any integer $k \in \Z$
introduce the spaces $Y_k$, defined by the
 norms\footnote{The expressions $\phi\to \norm{\phi}{Y_k}$
are not faithfull norms, in the sense that may be zero, even for
some  $\phi\neq 0$. On the other hand, they satisfy all the other
norm requirements and $\phi\to (\sum_k \|\phi_k\|_{Y_k}^2)^{1/2}$
is a norm!}
$$
\norm{\phi}{Y_k}= 2^{-k/2} \sum_m 2^{m/2} \norm{\phi_k}{L^2_t
L^2(|x|\sim 2^m)}.
$$
\\
Now we can define the Banach spaces $Y$ as a closure of the
functions
$$\phi(t,x) \in C_0^\infty( \R \times ( \R^n \setminus
 0))$$ with respect to the norm
\begin{equation}\label{eq.defY}
   \norm{\phi}{Y}:=\left(\suml_k \norm{\phi}{Y_k}^2\right)^{1/2}.
\end{equation}

Its dual space $Y^\prime$ consists of tempered distributions
$S^\prime(\R \times \R^n),$ having finite norm
$$\norm{\phi}{Y^\prime}:=\left(\suml_k
\norm{\phi}{Y^\prime_k}^2\right)^{1/2},
$$
where
$$
\norm{\phi}{Y^\prime_k}= 2^{k/2} \supl_m 2^{-m/2}
\norm{\phi_k}{L^2_t L^2(|x|\sim 2^m)}.
$$

Then the smoothing estimates \eqref{eq.3214}, \eqref{eq.sm2323}
and \eqref{eq.sminh2323} read

\begin{equation}\label{eq.3214Y}
   \norm{ e^{ i t \De} f}{Y^\prime}\leq C\norm{f}{L^2} \ ,\  \ \left\|\intl e^{ i s  \De}
F(s,\cdot) ds\right\|_{L^2_x}\leq C  \norm{F}{Y} .
\end{equation}
\begin{eqnarray}\label{eq.sminh2323Y}
    \left\| \intl_0^t e^{ i (t-s) \De} F(s,\cdot) ds\right\|_{Y^\prime}\leq   C  \norm{F}{Y} .
\end{eqnarray}

Motivated by the Strichartz estimates and these Besov versions of
``local smoothing'' norms, introduce the spaces
$$ X = L^1_t L^2_x + L^2_t L^{2n/(n+2)}_x + Y$$ with norm
$$
\|F\|_X =  \inf_{F=F^{(1)}+F^{(2)}+ F^{(3)}} \norm{F^{(1)}}{L^1_t
L^2_x} + \norm{F^{(2)}}{L^2_t L^{2n/(n+2)}_x} + \norm{F^{(3)}}{Y}
.$$

The dual to $X$ space is $X^\prime$ and the norm in this space is
defined in similar way:
\begin{equation}\label{eq.defXpr}
    \norm{\phi}{X^\prime}:=\left(\suml_k
\norm{\phi}{X^\prime_k}^2\right)^{1/2},
\end{equation}
where
$$
\norm{\phi}{X^\prime_k}= \supl_{(q,r)- Str.} \norm{\phi_k}{L^q_t
L^r_x}+ 2^{k/2} \supl_m 2^{-m/2} \norm{\phi_k}{L^2_t L^2(|x|\sim
2^m)}.
$$

The main result of this work is

\begin{thm}\label{mainSTR1} If $n \geq 3,$ then one can find a positive
number $\varepsilon >0$ so that for any (vector)   potential
$A=A(t,x)$ satisfying
\begin{eqnarray}\label{asA}
& &  \|\nabla A\|_{L^\infty L^{n/2}}+
\supl_k (\suml_m 2^m \norm{A_{<k} }{L^\infty L^\infty(|x|\sim
2^m)}) \leq \ve,
\end{eqnarray}
 there exists $C>0,$ such that for any
$F(t,x) \in S( \mathbb{R} \times \mathbb{R}^n)$ we have  the
estimate
\begin{eqnarray*}
 & &
  \  \left\| \  \intl_{t-s > 0} e^{ i
(t-s) \De_A} F(s,\cdot) ds \right\|_{X^\prime} \leq  C \| F\|_{X}
.
\end{eqnarray*}
In particular, the solutions to \eqref{eq.lwa} satisfy the
smoothing - Strichartz estimate
\begin{equation}\label{eq.mainsmstr}
\norm{u}{X^\prime} \lesssim \norm{f}{L^2}+ \norm{F}{X}.
\end{equation}
\end{thm}

\begin{remark}
The estimate \eqref{eq.mainsmstr} implies various interesting
inequalities. For example we have the classical Strichartz
estimate
$$
\supl_{(q,r)- Str.}\norm{u}{L^qL^r} \lesssim
\norm{f}{L^2}+\inf_{F=F_1+F_2}\norm{F_1}{L^1_t L^{2}_x} +
\norm{F_2}{L^2_t L^{2n/(n+2)}_x}
$$
as well as the smoothing - Strichartz estimates
$$
\norm{u}{Y^\prime} \lesssim
\norm{f}{L^2}+\inf_{F=F_1+F_2}\norm{F_1}{L^1_t L^{2}_x} +
\norm{F_2}{L^2_t L^{2n/(n+2)}_x},
$$
$$
\supl_{(q,r)- Str.}\norm{u}{L^qL^r} \lesssim
\norm{f}{L^2}+\norm{F}{Y} .
$$

\end{remark}

The main idea to prove this Theorem is to apply appropriate scale
invariant estimate for the free Schr\"odinger equation involving
Strichartz and smoothing type norms.

Estimates of this type have been obtained earlier in \cite{RuVe93}
and \cite{RuVe94} with Strichartz type norms of the form $\| F\|_{
L^{2n/(n+2)}_x L^2_t}.$ Recently, we found (the authors are
grateful to Luis Vega for pointing them this recent work) similar
estimate in the work \cite{IoK2005} and this estimate has the form
\begin{equation}\label{eq.S.4IK}
 \| D_x^{1/2}\intl_{0}^t e^{ i
(t-s) \De} F(s,\cdot) ds \|_{ \ L^\infty_{x_1} L^2_{x^\prime}
L^2_{t}} \leq
    C \| F\|_{L^2_t L^{2n/(n+2)}_x}.
\end{equation}
On one hand, this estimate  can be used to derive the Strichartz
estimate for the perturbed Schr\"odinger equation provided its
(formally)  "dual" version
\begin{equation}\label{eq.S.4dual}
 \left\|\ \intl_{0}^t e^{ i
(t-s) \De} F(s,\cdot) ds \right\|_{L^2_t L^{2n/(n-2)}_x } \leq
    C \|  D_x^{-1/2} F\|_{L^1_{x_1}
L^2_{x^\prime} L^2_{t}}.
\end{equation}
is verified. We apply \eqref{eq.S.4dual} and show that
\eqref{eq.mainsmstr} is satisfied for the free Schr\"odinger
equation. Once \eqref{eq.mainsmstr} is established for the free
case we show that these estimates are stable under small magnetic
perturbations satisfying \eqref{asA}. We conclude this section
with an interesting corollary, regarding the spectrum of $ \De_A$.
\begin{corollary}
Let $n\geq 3$ and
$A=A(x):\rn\to \rn$ be a real-valued vector potential, such that the
smallness conditions of Theorem \ref{mainSTR1} are met.
Let also $V=V(x):\rn\to \rone$, with $\norm{V}{L^{n/2}}<<1$.
  Then the spectrum
of $-\De_A+V$ does not contain eigenvalues.
\end{corollary}
This is a standard corollary of the Strichartz estimates in high dimensions.
Note that the requirement $n\geq 3$ is necessary,
and in fact such result fails in dimensions one and two.
\begin{proof}
Assume that there is an eigenvalue $\la$ with eigenvector $f$
for $-\De_A+V$. Then $u(t, \cdot)=
e^{i \la t} f$ is a solution to the Schr\"odinger equation
$$
u_t= i (-\De_A+V) u .
$$
It follows from the Strichartz estimates of Theorem \ref{mainSTR1}
\begin{eqnarray*}
& & \norm{u}{L^2(0,T) L^{2n/(n-2)}}\leq C( \norm{f}{L^2}+
\norm{V u}{L^2(0,T) L^{2n/(n+2)}})\leq \\
& & \leq  C( \norm{f}{L^2}+
\norm{V}{L^{n/2}} \norm{u}{L^2 L^{2n/(n-2)}}) \leq \\
& & C \norm{f}{L^2}+ C \ve \norm{u}{L^2 (0,T) L^{2n/(n-2)}} .
\end{eqnarray*}
Clearly, if $C\ve<1/2$, we have that
$$
\norm{u}{L^2(0,T) L^{2n/(n-2)}}\lesssim   \norm{f}{L^2},
$$
for every $T>0$, which is a impossible, since \\
$\norm{u}{L^2(0,T) L^{2n/(n-2)}}\geq C T^{1/2} \norm{f}{L^{2n/(n-2)}}$.
\end{proof}
In dimension two,
one may consider the Aharonov-Bohm type vector potentials
(i.e. of the form
$A(r, \theta)= g(r)\psi(\theta) (\sin(\theta), -\cos(\theta))$), for which
$-\De_A+V$ is unitarily equivalent to $-\De+V$, \cite{Balinsky}.
For the Schr\"odinger operators $-\De+V$ however, it is
well-known that eigenvalues may exists
for arbitrarily small (and smooth compactly)
 potentials $V$. This is due to B. Simon \cite{Simon1},
see also
\cite[p.274, Theorem XIII.80]{reedsimon4}.
Therefore, such result must  fail in
dimensions two.

Similar examples must be easier to construct in
dimension one, but as of now we are unaware of such  reference.\\
\\
{\it Acknowledgement:}
The authors are grateful to Herbert Koch and Nickolay Tzvetkov for
interesting discussions.

\newpage
\section{Preliminaries}
\subsection{Fourier transform and Littlewood-Paley projections}
\label{sec:2.1}
Define the Fourier transform  and its inverse by
\begin{eqnarray*}
& &
\hat{f}(\xi)= \intl_{\rn} f(x)e^{- 2\pi i x\cdot  \xi} dx\\
& &
 f(x)=\intl_{\rn} \hat{f}(\xi)e^{2\pi i x\cdot  \xi} d\xi
\end{eqnarray*}
Introduce a
positive, decreasing, smooth away from zero
 function $\chi:\rone_+ \to \rone_+ $, supported in
 $\{\xi:0\leq \xi \leq 2\}$ and $\chi(\xi)=1,$ $\forall 0 \leq \xi\leq 1.$
 Define $\vp(\xi)=\chi(\xi)-\chi(2\xi)$,
which is positive and
supported in the annulus $1/2\leq |\xi|\leq 2$.
We have that $\vp$ is smooth and
$\sum_{k\in \cz} \vp(2^{-k} \xi)=1$ for all $\xi\neq 0$.
In higher dimensions, we slightly abuse the notations and
denote a function with similar properties by the same name, i.e.
$\vp(\xi)=\vp(|\xi|)$, $\chi(x)=\chi(|x|)$ etc. Note that for $n>1$, $\chi(x):\rn\to\rone$ is a
smooth function even at zero.

The $k^{th}$
Littlewood-Paley projection is defined as a multiplier
type operator by
$\widehat{P_k f}(\xi)=\vp(2^{-k}\xi) \hat{f}(\xi)$.
Note that the kernel of $P_k$ is integrable, smooth and
real valued for every $k$.
In particular, it is bounded on every $L^p: 1\leq p\leq \infty$ and it
commutes with differential operators. Another helpful observation is that
for the differential operator $D_x^s$ defined via
 the multiplier $|\xi|^s$, one has
$$
D_x^s P_k u=2^{ks} \tilde{P}_k u,
$$
where $\tilde{P}_k$ is given by the multipleir $\tilde{\vp}(2^{-k}\xi)$, where
$\tilde{\vp}(\xi)= \vp(\xi)|\xi|^s$. \\
We also consider
$P_{<k}:=\sum_{l<k} P_l$, which essentially restricts the
Fourier transform to frequencies $\lesssim 2^k$.

Define also the function $\psi(\xi)=\chi(\xi/4)-\chi(4\xi)$. Note
that $\psi$ has similar support properties as $\vp$ and $\psi(\xi)
\vp(\xi)=\vp(\xi)$. Thus, we may also define the operators $Z_k$
by $\widehat{Z_k f}(\xi)=\psi(2^{-k} \xi) \hat{f}(\xi)$. By the
construction, $Z_k P_k=P_k$ and $Z_k=P_{k-2}+\ldots+
P_{k+1}$. \\
Recall a version of the
 Calder\'on commutator estimate (see for example Lemma
2.1 in  the work of Rodnianski and Tao, \cite{rodnianski}), which reads
$$
\norm{[P_k,f] g}{L^r}\leq C 2^{-k} \norm{\nabla f}{L^q}\norm{g}{L^p},
$$
whenever $1\leq r,p,q\leq \infty$ and $1/r=1/q+1/p$. \\
Also of interest will be the properties of products under the
action of $P_k$.
 Starting with the relations $$ P_k( f g) =
\sum_{\ell,m} P_k( f_\ell g_m),$$ $$ P_k( f_\ell g_m)  =
\intl_{\rn} \intl_{\rn} P_k(\xi) P_\ell(\xi-\eta)\hat{f}(\xi-\eta)
P_m(\eta)\hat{g}(\eta)e^{2\pi i x\cdot  \xi} d\xi d\eta,
$$ we exploit the property $ {\rm supp}\  P_k(\xi) \subseteq \{2^{k-1} \leq |\xi| \leq 2^{k+1}\} $ and see that
the sum can be restricted to the set $$ \{|\ell - m| \geq 2+N_0, \
| \max(\ell,m) - k| \leq 3\} \cup \{|\ell - m| \leq 1+N_0, \ k
\leq \max(\ell,m) +3 \},$$ where $N_0\geq 1$ is arbitrary number.
This domain can be enlarged slightly using the inequality  $
\max(\ell,m) \leq \ell +1 +N_0$ provided $ |\ell - m| \leq 1+N_0$.
So we can restrict the sum over the union of the following sets
(the first two are disjoint for $N_0\geq 5$, while the third one
can overlap with them)
 $$ \{m \leq k-N_0+1  , \
| \ell - k| \leq 3\} ,  \{\ell  \leq k-N_0+1, \ | m - k| \leq
3\}$$ and $$ \{|\ell - m| \leq 1+N_0, \ \ell \geq k-N_0-4 \}.$$

In conclusion, for any two (Schwartz ) functions $f, g$ we have
the pointwise estimate
\begin{eqnarray*}
|P_k( f g)(x)| &\leq & \suml_{l\geq k-N_0-4}\  \suml_{|m-\ell|\leq
1+N_0}|P_k ( f_l g _{m})(x)| + \\
&+& |P_k( f_{\leq k-N_0+1}\  g_{k-3\leq \cdot\leq  k+3})(x)| + \\
& &+ |P_k(  f_{k-3\leq \cdot\leq  k+3} \ g_{\leq k-N_0+1})(x)|
\end{eqnarray*}
Taking for determinacy $N_0=7,$ we get
\begin{eqnarray*}
& & | P_k( f g)(x)| \leq |f_{\leq k-6}(x) g_k(x)| + |[P_k, f _{\leq k-6}] g_{k-3\leq
\cdot \leq k+3}(x)| + \\
& & + |P_k(  f_{k-3\leq \cdot\leq  k+3} \
g_{\leq k-6})(x)| +\suml_{l\geq k-11}\  \suml_{|m-\ell|\leq 8}|P_k
( f_l g _{m})(x)|
\end{eqnarray*}
In particular, we  need an appropriate (product like!)
expression for \\
 $P_k (A\nabla u)$. The main term is clearly
when $\nabla u$ is in high frequency mode, while $\vec{A}$ is low
frequency. More precisely, according to our considerations above,
$$
P_k (A\nabla u)= A_{\leq k-6}\nabla u_k +E^k,
$$
where $E^k(x)$ satisfies the pointwise estimate
\begin{eqnarray}
\label{eq:1200} |E^k(x)| \leq   |[P_k, A_{\leq k-6}] \nabla
u_{k-3\leq \cdot k+3}(x)|  + \\ \suml_{l\geq k-11}\
\suml_{|m-\ell|\leq 8}|P_k ( A_l \cdot\nabla u_{m})(x)| +
\nonumber   | P_k ( A_{k-3\leq\cdot\leq k+3}\cdot \nabla u_{\leq
k-6})(x)|
\end{eqnarray}

Note that in terms of $L^p$ behavior and Littlewood-Paley theory,
one treats these error terms as if they were in the form
$(\partial_x A ) u$.

\subsection{Besov spaces versions of the ``local smoothing space''}

 The space $Y$ was introduced  as the closure of
$S(\R \times \R^n)$ with respect to the norm in \eqref{eq.defY},
where
\begin{equation}\label{eq.dellocsm}
    \norm{\phi}{Y_k}= 2^{k/2} \sum_m 2^{-m/2} \norm{P_k\phi}{L^2_t
L^2(|x|\sim 2^m)}.
\end{equation}
We can  replace $\|F\|_{L^2(|x|\sim 2^m)}$ by the comparable
expression  $ \|\vp(2^{-m}\cdot) F\|_{L^2}$.
 This will be done frequently (and without much discussions)
  in the sequel in order to make use of the Plancherel's theorem, which is of course valid only in the global $L^2$ space.
  We mention also that the norm $\|\phi\|_Y$ is scale invariant
  for rescale factors any diadic number.

 We show
that the  the ``local smoothing space'' defined as a closure of
Schwartz functions $\phi$ with respect to the ``local smoothing
norms'' $$\sum_m 2^{m/2} \norm{D_x^{-1/2} \phi(t,x)}{L^2_t
L^2(|x|\sim 2^m)}$$ can be embedded in  $Y.$
\begin{lem}
\label{le:4prima} There is a constant $C=C(n)$, so that for every
Schwartz function $\phi$ we have
\begin{equation}
\label{eq:09prima}\norm{\phi}{Y} \leq C  \sum_m 2^{m/2}
\norm{D_x^{-1/2} \phi(t,x)}{L^2_t L^2(|x|\sim 2^m)}.
\end{equation}
\end{lem}
\begin{proof} Taking into account the definition of the space $Y$, it is sufficient to establish the estimate
\begin{equation}\nonumber
\norm{\phi}{Y_k} \leq C  \sum_m 2^{m/2} \norm{D_x^{-1/2}
\phi_k(t,x)}{L^2_t L^2(|x|\sim 2^m)}
\end{equation}
for any integer $k.$ Using the scale invariance of the estimate we
see that we lose no generality taking $k=0.$ Thus, we have to
verify the estimate
$$  \sum_m 2^{m/2} \norm{ \vp(2^{-m}\cdot)P_0\phi}{L^2_t
L^2_x} \leq C \sum_m 2^{m/2} \norm{\vp(2^{-m}\cdot) D_x^{-1/2}
\phi_0(t,x)}{L^2_t L^2_x}.$$ Since
$$ P_0 \phi  = \sum_{|k|\leq 2}P_0 D_x^{1/2} D_x^{-1/2}P_k\phi =
\sum_{|k|\leq 2} \sum_{\ell \in \Z} \widetilde{P}_0
\vp(2^{-\ell}\cdot) D_x^{-1/2}P_k\phi,
$$
we can apply the triangle inequality, and reduce the proof to the
following estimate
$$  \sum_m 2^{m/2} \sum_{\ell \in \Z} \norm{ \vp(2^{-m}\cdot)\widetilde{P}_0
\vp(2^{-\ell}\cdot) D_x^{-1/2}P_k\phi}{L^2_t L^2_x} \leq  $$
$$ \leq C
\sum_\ell 2^{\ell/2} \norm{\vp(2^{-\ell}\cdot) D_x^{-1/2}
\phi_k(t,x)}{L^2_t L^2_x},$$ where $k \in \Z, |k| \leq 2.$ This
estimate follows easily from
\begin{equation}\label{Calstom1}
\norm{ \vp(2^{-m}\cdot)\widetilde{P}_0 \vp(2^{-\ell}\cdot) f}{
L^2_x} \leq C  \norm{f}{ L^2_x},
\end{equation}
\begin{equation}\label{Calstom2}
\norm{ \vp(2^{-m}\cdot)\widetilde{P}_0 \vp(2^{-\ell}\cdot) f}{
L^2_x} \leq C  2^{-m} \norm{f}{ L^2_x}, \ m \geq \ell+2
\end{equation}
 and the obvious observation that
$$ \sum_{m\leq \ell+1} 2^{m/2} +  \sum_{m\geq \ell +2} 2^{m/2}2^{-m} \lesssim 2^{\ell/2}.   $$
The estimate \eqref{Calstom1} is obvious, while the proof of
\eqref{Calstom2} follows from
$$ \vp(2^{-m}\cdot)\widetilde{P}_0 \vp(2^{-\ell}\cdot) f = [\vp(2^{-m}\cdot),\widetilde{P}_0] \vp(2^{-\ell}\cdot) f  ,
\ m \geq \ell +2 $$ and the Calder\'on estimate
$$ \norm{[\vp(2^{-m}\cdot),\widetilde{P}_0] g}{ L^2_x} \leq C 2^{-m} \norm{g}{
L^2_x}.$$ This completes the proof of the Lemma.
\end{proof}
\begin{remark} Note that the argument in the proof of this lemma
implies also the estimates
\begin{equation}\label{eq.bound}
\norm{P(D)f_k}Y \lesssim \norm{f_k}Y = \norm{f}{Y_k} , \ \forall k
\in \Z
\end{equation}
for any pseudodiferential operator with symbol $P(\xi) \in
C_0^\infty(\R^n).$
\end{remark}

We have also the  estimate (dual to \eqref{eq:09prima})
\begin{lem}
\label{le:4} There is a constant $C=C(n)$, so that for every
Schwartz function $\phi \in Y^\prime,$ we have
\begin{equation}
\label{eq:09mod} \supl_m 2^{-m/2} \norm{D_x^{1/2} \phi(t,x)}{L^2_t
L^2(|x|\sim 2^m)} \leq C_n \norm{\phi}{Y^\prime}.
\end{equation}
\end{lem}

\begin{remark} Some generalizations of the previous two Lemmas can
be seen in Theorem 1.6 and Theorem 1.7 in \cite{GeTa}.
\end{remark}

\section{Estimates for the bilinear form $Q(F, G)$}

The sesquilinear form
\begin{eqnarray*}
& & Q (F, G)=\int\int_{t>s}\langle e^{i(t-s)\Delta}F(s),
G(t)\rangle_{L^2(\R^n)} ds dt
\end{eqnarray*}
was used in \cite{KT} to derive Strichartz estimates (with
endpoint) and this estimates can be expressed in terms of $Q$
\begin{equation}
\label{eq:n2}
|Q(F,G)|\leq C\norm{F}{L^{q_1'}_t L^{r_1'}_x}  \norm{G}{L^{q_2'}_t L^{r_2'}_x},
\end{equation}
for all Strichartz pairs $(q_1, r_1)$, $(q_2, r_2)$.

We have the following estimate that can be
obtained by  applying Lemma 3 from the work of Ionescu-Kenig
\cite{IoK2005}.

\begin{thm}\label{KTapprpr}
There exists a constant $C=C(n)$ so that for any integer $k$, any
$F(t,x) \in S( \mathbb{R} \times \mathbb{R}^n)$ and  $G(t,x) \in
S( \mathbb{R} \times \mathbb{R}^n)$
\begin{equation}\label{T.25}
 |Q (F_k, G_k)|\leq C \left(\sum_{m \in
\mathbb{Z}}2^{m/2}2^{-k/2} \| \vp(2^{-m}\cdot)
 F_{k}\|_{L^2_t L^2_x} \right)\| G_k\|_{L^{2}_tL^{2n/(n+2)}_x}.
\end{equation}
\end{thm}

We have also the following energy-smoothing estimate.

\begin{thm}\label{KTensmqf}
There exists a constant $C=C(n)$ so that for any integer $k$, any
$F(t,x) \in S( \mathbb{R} \times \mathbb{R}^n)$ and  $G(t,x) \in
S( \mathbb{R} \times \mathbb{R}^n)$
\begin{equation}\label{T.25ensm}
 |Q (F_k, G_k)|\leq C \left(\sum_{m \in
\mathbb{Z}}2^{m/2}2^{-k/2} \| \vp(2^{-m}\cdot)
 F_{k}\|_{L^2_t L^2_x} \right)\| G_k\|_{L^{1}_tL^{2}_x}.
\end{equation}
\end{thm}

Before proving these  Theorems,  we verify some of the smoothing
estimates used in this work.

\subsection{Estimates in the local smoothing space}

For $n=1$ we have the following smoothing estimates (see Konig,
Ponce, Vega \cite{KPV91})
\begin{equation}\label{eq.S.1}
    2^{k/2} \|e^{-it\Delta} f_k \|_{L^\infty_x L^2_{\gamma}} \leq
    C \|f_k\|_{L^2},
\end{equation}
\begin{equation}\label{eq.S.1inh}
    2^{k/2} \left\|\int_{s<t} \ \ e^{-i(t-s)\Delta} F_k(s) ds \right\|_{L^\infty_x L^2_{t}} \leq
    C \|F_k\|_{L^1_x L^2_{t}},
\end{equation}
as well as
\begin{equation}\label{eq.S.2}
    2^{k/2} \left\| \int_\gamma \ \ e^{-it\Delta} F_k(t) dt \ \right\|_{L^2} \leq
    C \|F_k\|_{L^1_x L^2_{\gamma}}
\end{equation}
for any interval  $\gamma \subseteq \R_t.$ Here $C>0$ is a
constant independent of $f,F,\gamma.$

For $n >1$ we may assume
\begin{equation}\label{H1}
    {\rm supp} \ \ \widehat{f}(\xi) \subseteq \{ |\xi^\prime| \leq \xi_1/10, \ \ \xi^\prime = (\xi_2, \cdots, \xi_n)
    \}.
\end{equation}

Then we have the representation
\begin{eqnarray}\label{eqS3}
    \left(  e^{-it\Delta} f \right)(x_1, x^\prime) = \\ \nonumber = c
    \int_{\R^{n-1}} \int_{\R^{n-1}}  e^{it|\xi^\prime|^2 +
    i(x^\prime-y^\prime)\xi^\prime} \left( e^{-it\Delta_1} f\right)( x_1,
    y^\prime)\ d\xi^\prime d y^\prime,
\end{eqnarray}
where $\Delta_1 = \partial^2_{x_1}.$

This representation and one dimensional  estimates \eqref{eq.S.1},
\eqref{eq.S.1inh} and \eqref{eq.S.2}  lead to the following.
\begin{lem}
\label{le:1smn}  There exists a constant $C$ depending only on the
dimension, so that for any $f \in S(\R^n), F \in S(\R \times
\R^n),$ satisfying \eqref{H1} and
\begin{equation}\label{H2}
    {\rm supp}_\xi \ \ \widehat{F}(t,\xi) \subseteq \{ |\xi^\prime| \leq \xi_1/10, \ \ \xi^\prime = (\xi_2, \cdots, \xi_n)
    \}
\end{equation}
we have
\begin{equation}\label{eq.S.4}
    2^{k/2} \|e^{-it\Delta} f_k \|_{\ L^\infty_{x_1} L^2_{x^\prime}  L^2_{\gamma}} \leq
    C \|f_k\|_{L^2_{x^\prime} L^2_{x_1} }
\end{equation}
\begin{equation}\label{eq.S.4inh}
    2^{k/2} \left\|\int_{s<t} \ \ e^{-i(t-s)\Delta} F_k(s) ds \right\|_{L^\infty_{x_1} L^2_{x^\prime,t}} \leq
    C \|F_k\|_{L^1_{x_1} L^2_{x^\prime,t}},
\end{equation}
and
\begin{equation}\label{eq.S.6}
    2^{k/2} \left\| \int_\gamma \ \ e^{-it\Delta} F_k(t) dt \ \right\|_{L^2_{x_1} L^2_{x^\prime}} \leq
    C \|F_k\|_{ L^1_{x_1} L^2_{x^\prime} L^2_{\gamma}}
\end{equation}
for any interval $ \gamma \subseteq \R_t.$
\end{lem}

\begin{proof}  To prove \eqref{eq.S.4} we
use \eqref{eqS3} and find
\begin{eqnarray}\label{eqS3a}
    \left( \widehat{ e^{-it\Delta} f} \right)(x_1, \xi^\prime) = \\ \nonumber = c
      e^{it|\xi^\prime|^2 } \left( e^{-it\Delta_1} \widehat{f}\right)( x_1,
    \xi^\prime)= c
    e^{-it\Delta_1} \left(  e^{it|\xi^\prime|^2 }   \widehat{f}(\cdot, \xi^\prime) \right)(
    x_1).
\end{eqnarray}
Note that
\begin{equation}\label{eq.S5}
    P_k(\xi_1) \sim P_k(\xi)
\end{equation}
for $\xi \in {\rm supp}_\xi \ \ \widehat{f}$ due to \eqref{H1}.
>From this observation, the one dimensional estimate \eqref{eq.S.1}
and the Plancherel identity we get \eqref{eq.S.4}, since
$$
\|\hat{f}(x_1,\xi^\prime) \|_{L^2_{\xi^\prime} L^\infty_{x_1}}
\geq \|\hat{f}(x_1,\xi^\prime) \|_{ L^\infty_{x_1}
L^2_{\xi^\prime} } = \|f(x_1,x^\prime) \|_{ L^\infty_{x_1}
L^2_{x^\prime} } .
$$In a similar way we prove \eqref{eq.S.4inh} and \eqref{eq.S.6}.
This completes the proof
of the Lemma.
\end{proof}

Applying the H\"older inequalities
$$
\|g\|_{L^1_{x_1}} \lesssim \sum_{m \in \mathbb{Z}}2^{m/2} \|
g\|_{L^2_{x_1}(|x|\sim 2^m)},\ \supl_{m \in \mathbb{Z}}2^{-m/2} \|
g\|_{L^2_{x_1}(|x|\sim 2^m)} \leq \|g\|_{L^\infty_{x_1}}, $$ we obtain
\begin{corollary} \label{corsmvs}  The smoothing estimates \eqref{eq.3214},
\eqref{eq.sm2323}, \eqref{eq.sminh2323} are satisfied.
\end{corollary}

 By  Corollary \ref{corsmvs} one gets
\begin{equation}
\label{eq:121}
\begin{array}{l}
|Q(F_k, G_k)|\leq C_n (\suml_m 2^{-k/2} 2^{m/2}
\norm{F_{k}}{L^2_t L^2(|x|\sim 2^m)}) \times \\
\times (\suml_m 2^{-k/2} 2^{m/2} \norm{G_{k}}{L^2_t L^2(|x|\sim
2^m)} )
\end{array}
\end{equation}

After this prepartion,  we  turn to

\subsection{Proof of Theorem \ref{KTapprpr}: Bilinear smoothing-Strichartz estimate}

 The estimate \eqref{T.25} is scale invariant and for
this we can take $k=0.$ We have the relation
$$ Q (F, G) = \int\int_{\R^2}\langle e^{i(t-s)\Delta}F(s),
G(t)\rangle_{L^2(\R^n)} ds dt - $$ $$- \int\int_{t<s}\langle
e^{i(t-s)\Delta}F(s), G(t)\rangle_{L^2(\R^n)} ds dt$$ For the form
$$ Q_0 (F, G) = \left\langle \int_{\R} ds\  e^{-is\Delta}F(s), \int_{\R}
dt\  e^{-it\Delta} G(t) \right\rangle_{L^2(\R^n)} $$ we can apply
the Cauchy inequality and via \eqref{eq:201} and \eqref{eq.sm2323}
we get
$$ |Q_0(F_0, G_0)|\leq C \left(\sum_{m \in
\mathbb{Z}}2^{m/2} \| \vp(2^{-m}\cdot)
 F_{0}\|_{L^2_t L^2_x} \right)\| G_0\|_{L^{2}_tL^{2n/(n+2)}_x}.$$
Hence it remains to evaluate the form
$$ Q^*(F, G) =  \int\int_{t<s}\langle
e^{i(t-s)\Delta}F(s), G(t)\rangle_{L^2(\R^n)}$$ and verify the
inequality
\begin{equation}\label{T.25star}
 |Q^*(F_0, G_0)|\leq C \left(\sum_{m \in
\mathbb{Z}}2^{m/2} \| \vp(2^{-m}\cdot)
 F_{0}\|_{L^2_t L^2_x} \right)\| G_0\|_{L^{2}_tL^{2n/(n+2)}_x}.
 \end{equation}
To prove \eqref{T.25star} it is sufficient
to consider $F$ with
\begin{equation}\label{H2cone}
    {\rm supp}_\xi \ \ \widehat{F}(t,\xi)
\subseteq \{ |\xi^\prime| \leq \xi_1/10, \ \
\xi^\prime = (\xi_2, \cdots, \xi_n)
    \}.
\end{equation}
Also, note that
$$ Q^*(F, G) = \int_{\R \times \R^n} F(s,y) \overline{u(s,y)} \
ds dy,$$ where $u$ is a solution to the free Schr\"odinger
equation $ i\partial_t u + \Delta u = G$ having initial data
identically $0.$

With \eqref{H2cone} in mind, apply Lemma 3 in
Ionescu-Kenig \cite{IoK2005}. We get
\begin{equation}\label{eq.IK}
\| D^{1/2}_{x_1}
 u\|_{L^\infty_{x_1} L^2_{x^\prime,t}} \lesssim \|
 G\|_{L^{2}_tL^{2n/(n+2)}_x}
\end{equation}
Here and below we use the notations $x=(x_1,x^\prime),
x^\prime=(x_2,\cdots,x_n).$
 So we have
\begin{equation}\label{T.25pIK}
 |Q^*(F, G)|\leq C \left( \| D^{-1/2}_{x_1}
 F\|_{L^1_{x_1} L^2_{x^\prime,t}} \right)\| G\|_{L^{2}_tL^{2n/(n+2)}_x}.
\end{equation}
Thus, we need to  establish the inequality
$$ \| D^{-1/2}_{x_1}
 F_0\|_{L^1_{x_1} L^2_{x^\prime,t}} \lesssim \|F_0\|_{Y_0} = \sum_{m \in
\mathbb{Z}}2^{m/2} \| \vp(2^{-m}\cdot)
 F_{0}\|_{L^2_t L^2_x}. $$
 For the purpose it is
sufficient to apply \eqref{eq.bound}, the H\"older inequality
$$
\|g\|_{L^1_{x_1}} \lesssim \sum_{m \in \mathbb{Z}}2^{m/2} \|
g\|_{L^2_{x_1}},$$ and note that
$$ D^{-1/2}_{x_1}
 F_0 = P(D) F_0, $$
 for some $P(\xi) \in C_0^\infty(\R^n)$ due to our assumption
 \eqref{H2cone}. This completes the proof of the Theorem.

\subsection{Proof of Theorem \ref{KTensmqf}: bilinear energy --
smoothing estimate}

The proof follows the same line of the proof of Theorem
\ref{KTapprpr} with the following changement: in the place of
Ionescu-Kenig inequality \eqref{eq.IK} we use
\begin{equation}\label{eq.ensmooth23}
   \supl_t \left\|\int_0^t e^{ i (t-s)  \De}
F_k(s,\cdot) ds\right\|_{L^2_x}\leq C \left( \sum_{m \in \Z}
2^{m/2} 2^{-k/2} \norm{  F_k}{L^2_t L^2(|x|\sim 2^{m})} \right).
\end{equation}
This estimate is trivial, since by  the $L^2$ energy conservation,
 the left-hand  side of this inequality is equal to
$$\supl_t \left\|\int_0^t e^{ -i s  \De}
F_k(s,\cdot) ds\right\|_{L^2_x}$$  and applying the estimate
\eqref{eq.sm2323}, we can finish the proof as before.

\section{Proof of Theorem \ref{mainSTR1}}

 We start by some reductions of the problem. First,
note that \eqref{eq.lwa} is in the form
\begin{equation}\label{eq.lwa1}
\begin{cases}
\partial_t u - i\Delta u+2 A\nabla u = \tilde{F}(t,x) \\
u(0,x)=f(x),
\end{cases}
\end{equation}
where $\tilde{F}=F - div(A) u - i (\sum_j A_j^2) u$. We claim that suffices to prove
\begin{equation}
\label{eq:17} \norm{u}{X^\prime}\leq
C_n(\norm{f}{L^2}+\|\tilde{F}\|_{X}) ,
\end{equation}
for the solutions of \eqref{eq.lwa1}. Indeed, assuming the validity of \eqref{eq:17}
 and since by our assumptions and Sobolev embedding
$ \norm{\nabla A}{L^\infty_t L^{n/2}_x}+\norm{A}{L^\infty_t L^{n}_x}\leq
C \norm{\nabla A}{L^\infty_t L^{n/2}_x}\leq C \ve$. We have
\begin{eqnarray*}
& &\norm{u}{X^\prime}\leq C\norm{f}{L^2}+C  \|\tilde{F}\|_{X}\leq \\
& &\leq C\norm{f}{L^2}+C \norm{F}{X}+ C(\norm{\nabla A}{L^\infty_t
L^{n/2}_x}+\norm{A}{L^\infty_t L^{n}_x}^2)
\norm{u}{L^2L^{2n/(n-2)}}\\
& &\leq C_n\norm{f}{L^2}+C_n \norm{F}{X}+
C_n \ve \norm{u}{L^2L^{2n/(n-2)}}\leq \\
& & \leq  C_n\norm{f}{L^2}+C_n \norm{F}{X}+ C_n \ve
\norm{u}{X^\prime}.
\end{eqnarray*}
It follows that
$$
\norm{u}{X^\prime}\leq C \norm{f}{L^2}+C \norm{F}{X},
$$
as claimed, as long as $\ve: C_n\ve<1/2$.

Thus, we concentrate on showing \eqref{eq:17} for the solutions of
\eqref{eq.lwa1}, where we denote the right hand side by $F$ again.

Next, we take a Littlewood-Paley projection of \eqref{eq.lwa1}. We get
$$
\partial_t u_k - i\Delta u_k =
F_k-2 A_{<k-6}\nabla u_{k} - 2 E^k:=H_k,
$$
where $E^k$ is the error term $E^k=P_k(A\nabla u)-
A_{\leq k-6}\nabla u_k$ given by \eqref{eq:1200}.

We will show that the solution to
$\partial_t u_k - i\Delta u_k=H_k$ with
initial data $u_k(0,x)=f_k$, satisfies the estimate
\begin{equation}
\label{eq:18}
\begin{array}{c}
\norm{u_k}{X^\prime}\leq C\norm{f_k}{L^2}+  C \|H_k\|_{X}
\end{array}
\end{equation}

We will
 show first how  \eqref{eq:18} implies Theorem \ref{mainSTR1} and then
 we proceed to show \eqref{eq:18}.
\subsection{\eqref{eq:18} implies Theorem \ref{mainSTR1}}
Apply \eqref{eq:18} to $u_k$. We have
\begin{equation}
\label{eq:35}
\begin{array}{l}
\norm{u_k}{X^\prime}\leq C\norm{f_k}{L^2}+C(\norm{F_k}{X}+
\norm{E^k}{L^2 L^{2n/(n+2)}}+ \\
+ C \suml_m 2^{m/2} 2^{-k/2} \norm{A_{<k-6}\nabla
u_{k}}{L^2L^2(|x|\sim 2^m)})
\end{array}
\end{equation}
We will need the following estimates.
\begin{eqnarray}
\label{eq:30}
& & (\suml_k \|E^k\|_{L^2L^{2n/(n+2)}}^2)^{1/2}\leq C_n \ve
(\suml_k \norm{u}{L^2 L^{2n/(n-2)}}^2)^{1/2} \leq
C_n \ve \norm{u}{X^\prime} \\
\label{eq:31} & & \suml_m 2^{m/2} 2^{-k/2} \norm{A_{<k-6}\nabla
u_{k}}{L^2L^2(|x|\sim 2^m)}\leq C_n\ve
\norm{u_k}{X^\prime_k}.
\end{eqnarray}
Let us show first how based on \eqref{eq:30} and  \eqref{eq:31},
we finish the proof of Theorem \ref{mainSTR1}. Plugging in these
estimates in \eqref{eq:35}, using the definition \eqref{eq.defXpr}
of $X^\prime$  and square summing in $k$ yields
\begin{eqnarray*}
& & \norm{u}{X^\prime}=\left(\suml_k
\norm{u_k}{X^\prime_k}^2\right)^{1/2}\leq C_n( \norm{f}{L^2}+
\norm{F}{X})+ C_n\ve \norm{u}{X^\prime},
\end{eqnarray*}
whence since $\ve: C_n\ve<1/2$,
$$
 \norm{u}{X^\prime}\leq C_n ( \norm{f}{L^2}+ \norm{F}{X}).
$$
Thus, for this section, remains to see \eqref{eq:30} and \eqref{eq:31}.
\subsubsection{Proof of \eqref{eq:31}} Let
$\widetilde{k}$ be integer with $|k-\widetilde{k}| \leq 3.$ We have
\begin{eqnarray*}
& & \suml_m 2^{m/2} 2^{-k/2} \norm{A_{<k-6}\nabla
u_{\widetilde{k}}}{L^2L^2(|x|\sim 2^m)}\lesssim \\
& &\lesssim
\left(\suml_m 2^{m} \norm{A_{<k-5}}{L^\infty L^\infty(|x|\sim 2^m)}
\right) \supl_m 2^{-m/2} 2^{-k/2} \norm{\nabla u_{\widetilde{k}}}{L^2L^2(|x|\sim 2^m)}\leq \\
& & \leq C_n  \ve \supl_m 2^{-m/2} 2^{-k/2} \norm{\nabla
u_{\widetilde{k}}}{L^2L^2(|x|\sim 2^m)}.
\end{eqnarray*}
This last expression is very similar to
$\norm{u_{\widetilde{k}}}{X^\prime}$.
We will show that it is controlled by it,
which of course is enough to establish \eqref{eq:31}.  \\
Fix an $m$. Then
\begin{eqnarray*}
& &  2^{-m/2} 2^{-k/2} \norm{\nabla
u_{\widetilde{k}}}{L^2L^2(|x|\sim 2^m)}\lesssim 2^{-m/2} 2^{k/2}
\norm{\vp(2^{-m}\cdot) Q_{\widetilde{k}}
u_{\widetilde{k}}}{L^2L^2},
\end{eqnarray*}
where $Q_k$ acts as a (vector) multiplier $\psi(2^{-k}\xi) 2^{-k}\xi$.
We have by the Calder\'on commutator estimate\footnote{We are
using the particular
form
 $ \norm{[Q_k,\vp(2^{-m}\cdot)]u_k}{L^2L^2}\lesssim 2^{-k} 2^{-m}
 \norm{(\nabla \vp)(2^{-m}\cdot)}{L^{2n/3}}
\norm{u_k}{L^{2n/(n-3)}}=2^{-k} 2^{m/2}
 \norm{u_k}{L^{2n/(n-3)}}$}
 and the Bernstein inequality
 \begin{eqnarray*}
& & 2^{-m/2} 2^{k/2} \norm{\vp(2^{-m}\cdot) Q_k u_k}{L^2L^2}\leq
2^{-m/2} 2^{k/2} \norm{Q_k(\vp(2^{-m}\cdot) u_k)}{L^2L^2}+ \\
& & +
2^{-m/2} 2^{k/2} \norm{[Q_k,\vp(2^{-m}\cdot)]u_k}{L^2L^2}\lesssim
2^{-m/2} 2^{k/2} \norm{\vp(2^{-m}\cdot) u_k}{L^2L^2}+\\
& &+ 2^{-k/2} \norm{u_k}{L^2_t L^{2n/(n-3)}}\lesssim
 2^{-m/2} 2^{k/2} \norm{\vp(2^{-m}\cdot) u_k}{L^2L^2}+ \norm{u_k}{L^2_t L^{2n/(n-2)}} \\
 & &\leq C_n \norm{u_k}{X_k^\prime}.
\end{eqnarray*}
\subsubsection{Proof of \eqref{eq:30}}
We treat $E^k$ on a term-by-term basis in \eqref{eq:1200}. For the first term, by Calder\'on commutators,
\begin{eqnarray*}
& & (\suml_k \norm{[P_k, A_{< k-6}] \nabla u_{k}
}{L^2L^{2n/(n+2)}}^2)^{1/2}
\lesssim \\
& & \lesssim (\suml_k \norm{\nabla A_{< k-6}}{L^\infty L^{n/2}}^2
\norm{u_{k-3\leq \cdot k+3}}{L^2L^{2n/(n-2)}}^2)^{1/2}\lesssim \\
& & \lesssim \supl_k
\norm{\nabla A_{< k-6}}{L^\infty L^{n/2}}
(\suml_k \norm{u_{k-3\leq \cdot k+3}}{L^2L^{2n/(n-2)}}^2)^{1/2}\lesssim \\
& &\lesssim
\norm{\nabla A}{L^\infty L^{n/2}}
\norm{u}{X^\prime}.
\end{eqnarray*}
 For the second
term, we have by standard Littlewood-Paley theory
\begin{eqnarray*}
& &
(\suml_k \norm{P_k G}{L^2 L^{2n/(n+2)}}^2)^{1/2}\lesssim \norm{G}{L^2 L^{2n/(n+2)}}, \\
& &
\|(\suml_l |g_l|^2)^{1/2}\|_{L^p}\sim \norm{g}{L^p} \qq \textup{for all}\q 1<p<\infty.
\end{eqnarray*}
whence with $m,\ell \in \Z$ with $|m-\ell|\leq 8$ we have
\begin{eqnarray*}
& & (\suml_k \|P_k(\suml_{|\ell-m|\leq 8} A_\ell \cdot\nabla u _{m})\|_{L^2L^{2n/(n+2)}}^2)^{1/2} \sim \\
& & \sim
\|\suml_{|\ell-m|\leq 8} A_\ell \cdot\nabla u _{m}\|_{L^2L^{2n/(n+2)}}\lesssim \\
& & \lesssim \|(\suml_\ell 2^{2l} |A_\ell|^2)^{1/2}\|_{L^\infty
L^{n/2}} \|(\suml_m |\tilde{P}_{m}
u|^2)^{1/2}\|_{L^2 L^{2n/(n-2)}} \sim \\
& &\sim \norm{\nabla A}{L^\infty L^{n/2}} \norm{u}{L^2
L^{2n/(n-2)}}\lesssim \ve \norm{u}{X^\prime}.
\end{eqnarray*}

For the third term in \eqref{eq:1200}, observe that since for all
$1\leq p\leq 2$,
$$
(\suml_k \|G^k\|_{L^p}^2)^{1/2}\leq C_n \|(\suml_k |G^k|^2)^{1/2}\|_{L^p}
$$
we can estimate by
\begin{eqnarray*}
& & (\suml_k \norm{P_k ( A_{k-3\leq\cdot\leq k+3}\cdot \nabla u_{< k-5})}{L^2L^{2n/(n+2)}}^2)^{1/2}\lesssim  \\
& & \lesssim (\suml_k 2^{2k}
\|A_{k-3\leq\cdot\leq k+3}  \tilde{P}_{< k-5} u \|_{L^2L^{2n/(n+2)}}^2)^{1/2} \lesssim \\
& & \lesssim\| (\suml_k 2^{2k}  |A_{k-3\leq\cdot\leq k+3}|^2 |\tilde{P}_{< k-5} u|^2)^{1/2}\|_{L^2L^{2n/(n+2)}} \lesssim \\
& & \lesssim \| (\suml_k 2^{2k}  |A_{k-3\leq\cdot\leq
k+3}|^2)^{1/2}\|_{L^\infty L^{n/2}}
 \|\supl_k|\tilde{P}_{< k-5} u|\|_{L^2L^{2n/(n-2)}}\lesssim \\
& & \lesssim \|\nabla A\|_{L^\infty L^{n/2}}
\norm{u}{L^2L^{2n/(n-2)}}\lesssim \ve \norm{u}{X^\prime},
\end{eqnarray*}
Here, we have used the pointwise estimate(see section 6.1, Chapter
I, \cite{Stei1})  $\sup_k|\tilde{P}_{< k-5} u|(x)\leq C M(u)(x)$,
where $M(u)$ is
the Hardy-Littlewood maximal function and therefore  \\
$$\|\sup_k|\tilde{P}_{< k-5} u|\|_{L^p}\leq C\|u\|_{L^p}$$ for all
$1<p<\infty$.

\subsection{The proof of \eqref{eq:18}}

 The nontrivial part of \eqref{eq:18} is
the case when $u_k$ is the solution to $\partial_t u_k - i\Delta
u_k=H_k$ with zero initial data $u_k(0,x)=0$. Then the fact that
the norm of $X_k$ has three components implies that the inequality
$$ \norm{u_k}{X^\prime}\leq  C \|H_k\|_{X}  $$
is equivalent to the following nine inequalities
\begin{equation}
\label{eq:100nac}  \norm{u_k}{L^2 L^{2n/(n-2)}}\leq C
\norm{H_k}{L^2_t L^{2n/(n+2)}}
\end{equation}
\begin{equation}
\label{eq:100nacen}  \norm{u_k}{L^2 L^{2n/(n-2)}}\leq C
\norm{H_k}{L^1_t L^{2}}
\end{equation}
\begin{equation}
\label{eq:100le}  \norm{u_k}{L^2 L^{2n/(n-2)}}\leq C \suml_m
2^{m/2} 2^{-k/2} \norm{H_k}{L^2_t L^{2}(|x|\sim 2^m)},
\end{equation}

\begin{equation}
\label{eq:100nacee}  \norm{u_k}{L^\infty L^{2}}\leq C
\norm{H_k}{L^2_t L^{2n/(n+2)}}
\end{equation}
\begin{equation}
\label{eq:100naceneee}  \norm{u_k}{L^\infty L^{2}}\leq C
\norm{H_k}{L^1_t L^{2}}
\end{equation}
\begin{equation}
\label{eq:100leeee}  \norm{u_k}{L^\infty L^{2}}\leq C \suml_m
2^{m/2} 2^{-k/2} \norm{H_k}{L^2_t L^{2}(|x|\sim 2^m)},
\end{equation}

\begin{equation}
\label{eq:100}
\begin{array}{c}
2^{k/2} \sup_m 2^{-m/2} \norm{u_k}{L^2_t L^{2}(|x|\sim 2^m)} \leq
C \norm{H_k}{L^2_t L^{2n/(n+2)}},
\end{array}
\end{equation}
\begin{equation}
\label{eq:100smen}
\begin{array}{c}
2^{k/2} \sup_m 2^{-m/2} \norm{u_k}{L^2_t L^{2}(|x|\sim 2^m)} \leq
C \norm{H_k}{L^1_t L^{2}},
\end{array}
\end{equation}
and
\begin{equation}
\label{eq:101}
\begin{array}{c}
2^{k/2} \sup_m 2^{-m/2} \norm{u_k}{L^2_t L^{2}(|x|\sim 2^m)}\leq \\
\leq C \suml_m 2^{m/2} 2^{-k/2} \norm{ H_k}{L^2_t L^{2}(|x|\sim
2^m)} .
\end{array}
\end{equation}

The estimates \eqref{eq:100nac}, \eqref{eq:100nacen},
\eqref{eq:100nacee} and \eqref{eq:100naceneee} are Strichartz
inequalities (see \eqref{eq:202} for general case).

The estimate \eqref{eq:101} is smoothing - smoothing estimate
established in Corollary \ref{corsmvs}( actually they follow from
the bilinear estimate \eqref{eq:121}) .

The estimates  \eqref{eq:100le}, \eqref{eq:100} are smoothing -
endpoint Strichartz inequalities following from bilinear estimate
of Theorem \ref{KTapprpr}.

Finally, the estimates  \eqref{eq:100leeee}, \eqref{eq:100smen}
are smoothing - energy inequalities following from bilinear
estimate of Theorem \ref{KTensmqf}.

The first inequality is the usual Strichartz estimate, while the second one is equivalent to \eqref{T.25}.

This completes the proof of the inequality \eqref{eq:18} and
Theorem \ref{mainSTR1}.

\end{document}